\magnification=\magstep1
\input pictex
\input amstex
\UseAMSsymbols
\NoBlackBoxes
       \font\kl=cmr8    
        
       \font\gr=cmbx10 scaled\magstep1 

\def\mod{\operatorname{mod}}
\def\Ext{\operatorname{Ext}}

\def\Hom{\operatorname{Hom}}

\def\arr#1#2{\arrow <2mm> [0.25,0.75] from #1 to #2}
\def\Arr#1#2{\arrow <2mm> [0.25,0.75] from #1 to #2}

\def\Rahmen#1%
   {$$\vbox{\hrule\hbox%
                  {\vrule%
                       \hskip0.5cm%
                            \vbox{\vskip0.3cm\relax%
                               \hbox{$\displaystyle{#1}$}%
                                  \vskip0.3cm}%
                       \hskip0.5cm%
                  \vrule}%
           \hrule}$$}           
\def\rahmen#1%
   {\vbox{\hrule\hbox%
                  {\vrule%
                       \hskip0.5cm%
                            \vbox{\vskip0.3cm\relax%
                               \hbox{$\displaystyle{#1}$}%
                                  \vskip0.3cm}%
                       \hskip0.5cm%
                  \vrule}%
           \hrule}}           

\centerline{\gr The Self-injective Cluster Tilted Algebras}
	\bigskip
\centerline{Claus Michael Ringel}
	\bigskip
{\narrower\narrower\narrower {\kl Abstract. We are going to determine all the self-injective
cluster tilted algebras. All are of finite
representation type and special biserial. There are two different
classes. The first class are the serial (or Nakayama) algebras with $\ssize n \ge 3$ 
simple modules and Loewy length $\ssize n-1$. The second class of algebras 
has an even number $\ssize 2m$ of simple modules, they are exhibited by quivers and relations 
at the end of this note.}
\par}
	\bigskip
Let $k$ be a field. The algebras we will deal with are finite-dimensional associative
$k$-algebras with 1. Given such an algebra $A$, the modules considered will usually
be (finite-dimensional)
left $A$-modules. The $k$-duality will be denoted by $D = \Hom(-,k).$
	
The class of cluster tilted algebras has been introduced by Buan, Marsh and Reiten [BMR].
According
to Zhu [Z] and Assem, Br\"ustle and Schiffler [ABS], they can be defined as the semi-direct extension
of $A$ by the $A$-$A$-bimodule $\Ext^2(D(A_A),{}_AA)$, where $A$ is a tilted algebra.
The aim of
this note is to single out those cluster tilted algebras which are self-injective. 
It is rather easy to see that self-injective cluster tilted algebras have to be of finite
representation type. We will see that all are special biserial. There are two different
classes. The first class are the serial (or Nakayama) algebras with $n\ge 3$ simple modules and Kupisch series $(n-1,\dots,n-1)$. The second class of algebras 
has an even number $2m$ of simple modules, they are exhibited by quivers and relations 
at the end of this note.
	\bigskip
{\bf First considerations.}
	\medskip
Let $H$ be a connected hereditary finite-dimensional $k$-algebra. Let $n$ be the number
of isomorphism classes of simple $H$-modules. We
denote by $\Cal C(H)$
the corresponding cluster category; by definition, this is the orbit category
of the derived category $D^b(\mod H)$ with respect to the functor $F = \tau_d^{-1}[1]$, where
$[1]$ is the shift functor and $\tau_d$ the Auslander-Reiten translation functor in the derived
category. We denote by $\tau_c$ the Auslander-Reiten translation functor in the cluster
category; both $\tau_d$ and $\tau_c$ are invertible. The Auslander-Reiten translation
in $\mod H$ itself will be denoted by $\tau_H$.

Let $T$ be a multiplicity-free tilting $H$-module, we may consider it as an object in the cluster
category, and we denote by $\widetilde A$ its endomorphism ring  in the cluster
category. This is a typical cluster tilted algebra. Write 
$T = \bigoplus_{i=1}^n T_i$ with $T_i$ indecomposable. 

The question which we consider here is: when is  $\widetilde A$ self-injective? 
Note that {\it $\widetilde A$ is self-injective if and only if 
$\tau_c^2T$ is isomorphic to $T$.}
Namely, the category $\mod \widetilde A$ is the factor category
of $\Cal C(H)$ modulo the ideal generated by the objects $\tau_cT_i$; in this factor category,
the objects $T_i$ are the indecomposable projective ones, the objects $\tau_c^2 T_i$ the
indecomposable injective ones. 
In order for $\widetilde A$ to be self-injective we just need
 that
any indecomposable injective $\widetilde A$-module (thus any $\tau_c^2T_i$) is projective
(thus isomorphic to some $T_j$).
	\bigskip
\vfill\eject

{\bf Exclusion.}
	\medskip

Assume that $T$ and $\tau^2_cT$ are isomorphic.
Now $T$ has only finitely many indecomposable direct summands (namely $n$); this means
that any object $T_i$ has to be $\tau_c$-periodic. If $H$ is wild, then $\Cal C(H)$ has no
$\tau_c$-periodic objects at all. Thus $H$ cannot be wild. If $H$ is tame, then the
$\tau_c$-periodic objects in $\Cal C(H)$ are actually $\tau_H$-periodic modules. But it is
well-known that the number of $\tau_H$-periodic indecomposable
direct summands of 
a tilting $H$-module is at most $n-2$.  (A stable tube of rank $r$ can have only $r-1$ indecomposable direct summands
which are pairwise $\Ext$-orthogonal, and the sum of these numbers $r-1$ over all the tubes
is $n-2$, see [DR].) 
This shows that {\it $H$ has to be representation-finite.}

Note that the $\tau_c$-orbit of any summand $T_i$ has to have an even number of elements,
so that half of these elements can be direct summands of $T$. 

We have to consider now the various cases $A_n, B_n,\dots,E_8, G_2$ in detail.
Let
$\Delta$ be the quiver of $H$ (we consider it as a valued quiver (see [DR]),
if we deal with one of the cases $B_n, C_n, F_4$
and $G_2$). Let $a$ be a vertex of $\Delta$. Denote by $P(a)$ the
indecomposable projective $H$-module corresponding to the vertex $a$.  
We are looking for natural numbers $t$ such that $\Hom(P(a),\tau_H^{-t}P(a)) \neq 0.$
If $\tau_H^{-t}P(a)$ is not injective, then we see that $\Ext^1(\tau_H^{-t-1}P(a),P(a))\neq 0.$
Now, if $t$ is odd, then $-t-1$ is even: We conclude that in this case neither
$P(a)$ nor any other element of the $\tau_c$-orbit of $P(a)$ can be a direct summand of $T$.
In order to decide whether $\Hom(P(a),\tau_H^{-t}P(a))$ is zero or not, one just has
to calculate the hammock function starting at $P(a)$, see [G] or [RV].

Let us assume now that  $T$ has an indecomposable direct summand $T_i$ which is in the
$\tau_H$-orbit of $P(a)$.

(1) {\it The vertex $a$ is a boundary vertex of $\Delta$.}
Assume for the contrary that $a$ is an interior vertex,
We have $\Hom_H(P(a),\tau^{-1}P(a)) \neq 0$. If $\Delta$ is different from $A_3$, then
$\tau^{-1}P(a)$ is not injective, thus we get a contradiction.
If $\Delta$ is of type $A_3$, then the orbit of $\tau_c$ corresponding to $a$
has length 3, thus is of odd length, again a contradiction.

(2) {\it The edge between $a$ and its (unique) neighbor has valuation $(1,1)$.}
Here again we see that otherwise $\Hom_H(P(a),\tau^{-1}P(a)) \neq 0$, and in this case 
$\tau^{-1}P(a)$ cannot be injective. This immediately excludes $G_2$. Also, together with (1)
it excludes $B_n$ and $C_n$: Namely, all the indecomposable summands of $T$ would have
to belong to a single $\tau_H$-orbit --- but since the simple $H$-modules do have two different
kinds of endomorphism rings, the same is true for any tilting $H$-module ([R], alternatively,
we also could argue that the orbit in question has not enough elements).

Assume that $\Delta$ is of type $D_n$ or $E_n$, and that $a = a_1, a_2,\dots,a_p$
is the minimal path from $a$ to the branching vertex $a_p$. If necessary, we write $p = p(a).$ 

(3) {\it The number $p$ is even.} (We have $\Hom_H(P(a),\tau^{-p}P(a)) \neq 0$, and 
$\tau^{-p}P(a)$ is not injective). 

(4) {\it If $\Delta = E_n$, then $p > 2.$} (Again, $\Hom_H(P(a),\tau^{-3}P(a)) \neq 0$, and 
$\tau^{-3}P(a)$ is not injective). 

The assertions (1), (3) and (4) exclude immediately the cases $E_6, E_8$, and they show that in 
case $E_7$ all the indecomposable direct summands of $T$ have to belong to the $\tau_c$-orbit
which contains the module $P(a)$ with $p(a) =  4.$ But any $\tau_c$-orbit of $\Cal C(H)$ with
$H$ of type $E_7$ is of length 10, thus  only 5 indecomposable direct summands of $T$
can belong to this orbit. Thus also $E_7$ is excluded.
Similarly, we see that 
in case $H$ is of type $D_n$, and $a(p) > 2$, then $n = a(p)+2$ has
to be even.

In order to exclude the case $F_4$, we only have to notice that 
 $\Hom(P(a),\tau^{-3}P(a)) \neq 0$
for any of the two boundary vertices.

We still have to consider the cases $A_n$. The two $\tau_H$-orbits of $H$ yield a single
$\tau_ c$-orbit containing $n+3$ elements. These are the only possible summands of $T$,
thus we must have $\frac12(n+3) = n$, therefore $n = 3.$ 

Besides $A_3$ also the cases $D_n$ remain. Assume we are in case $D_n$ and $p(a) > 2.$ The
$\tau_c$-orbit of $P(a)$ contains precisely $n$ elements, this shows that $n$ has to be 
even. For $n \ge 5$, there are two $\tau_H$-orbits $a$ with $a(p) = 2$; for $n = 4$, there
are three such $\tau_H$-orbits.
In case $n$ is odd, the two $\tau_H$-orbits with $a(p) = 2$
combine to form part of a single $\tau_c$-orbit with $2n$ elements. In case $n$ is even,
the $\tau_H$-orbits with $a(p) = 2$ 
are contained in two (or, for $n=4$, three) separate $\tau_c$-orbits, each having $n$ elements. 
All these $\tau_c$-orbits actually occur, as we are going to show next.
	\bigskip
{\bf The possible cases.}
	\medskip
{\bf 1. Case $A_3$.} The algebra $\widetilde A$ is serial with 3 simple modules and
Kupisch series $(2,2,2)$.
Here are the Auslander-Reiten quivers of $\Cal C(H)$ and of 
$\mod \widetilde A$. Note that in both cases the left hand boundary has to be 
identified with the right hand side using a twist in order to form a M\"obius strip. The summands $T_i$
are marked by a star $*$.

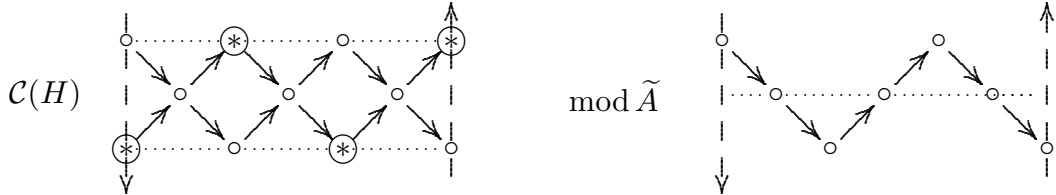
$$
\hbox{\beginpicture
\setcoordinatesystem units <.6cm,.6cm>
\put{\beginpicture
\multiput{$\circ$} at 0 2  1 1  2 0  3 1  4 2  5 1  6 0 /
\multiput{$*$} at 0 0  2 2  4 0  6 2 /
\multiput{$\bigcirc$} at 0 0  2 2  4 0  6 2 /
\arr{0.2 0.2}{0.8 0.8}
\arr{0.2 1.8}{0.8 1.2}
\arr{1.2 0.8}{1.8 0.2}
\arr{1.2 1.2}{1.8 1.8}

\arr{2.2 0.2}{2.8 0.8}
\arr{2.2 1.8}{2.8 1.2}
\arr{3.2 0.8}{3.8 0.2}
\arr{3.2 1.2}{3.8 1.8}

\arr{4.2 0.2}{4.8 0.8}
\arr{4.2 1.8}{4.8 1.2}
\arr{5.2 0.8}{5.8 0.2}
\arr{5.2 1.2}{5.8 1.8}
\setdots <1mm>
\plot 0.2 0  1.8 0 /
\plot 2.2 0  3.8 0 /
\plot 4.2 0  5.8 0 /

\plot 0.2 2  1.8 2 /
\plot 2.2 2  3.8 2 /
\plot 4.2 2  5.8 2 /
\setdashes <2mm>

\arr{0 2.5}{0 -.8}
\arr{6 -.5}{6 2.7}
\endpicture} at 0 0
\put{\beginpicture
\multiput{$\circ$} at 0 2  1 1  2 0  3 1  4 2  5 1  6 0 /
\arr{0.2 1.8}{0.8 1.2}
\arr{1.2 0.8}{1.8 0.2}

\arr{2.2 0.2}{2.8 0.8}
\arr{3.2 1.2}{3.8 1.8}

\arr{4.2 1.8}{4.8 1.2}
\arr{5.2 0.8}{5.8 0.2}
\setdots <1mm>
\plot 0.2 1  0.8 1 /
\plot 1.2 1  2.8 1 /
\plot 3.2 1  5.8 1 /

\plot 5.2 1  5.8 1 /
\setdashes <2mm>

\arr{0 2.5}{0 -.8}
\arr{6 -.5}{6 2.7}
\endpicture} at 11 0
\put{$\Cal C(H)$} at -4.5 0
\put{$\mod \widetilde A$} at 6 0
\endpicture}
$$
	\medskip
{\bf 2. The case $D_n$ with $n\ge 4$, using two short arms.} Let us present the
cases $n = 7$ and $n = 8$.

First, we present the case $n=7.$ In this case, the upper part (given by  the two short arms)
provides a M\"obius strip (the lower part given by the long arm is always a cylinder). In
general, for $n$ odd one, the short arm part of the Auslander-Reiten quiver is a M\"obius strip.
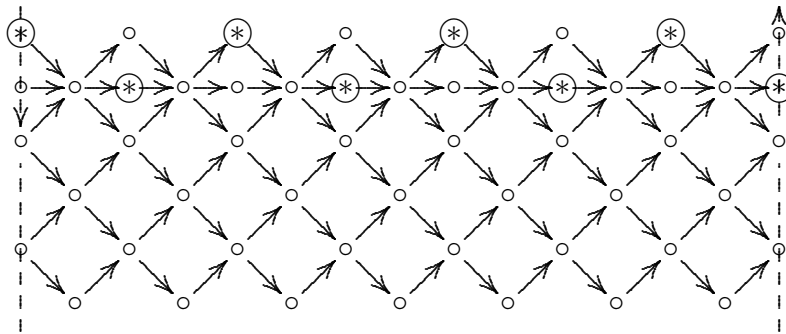
$$
\hbox{\beginpicture
\setcoordinatesystem units <.6cm,.6cm>
\multiput{$\circ$} at 
       0 2  0 4  0 5  
  1 1  1 3  1 5
       2 2  2 4       2 6
  3 1  3 3  3 5
       4 2  4 4  4 5  
  5 1  5 3  5 5
       6 2  6 4       6 6
  7 1  7 3  7 5
       8 2  8 4  8 5  
  9 1  9 3  9 5
        10 2  10 4        10 6
  11 1  11 3  11 5
        12 2  12 4  12 5  
  13 1  13 3  13 5
        14 2  14 4        14 6 /
\multiput{$*$} at 0 6  2 5  4 6  6 5  8 6  10 5  12 6  14 5 /
\multiput{$\bigcirc$} at 0 6  2 5  4 6  6 5  8 6  10 5  12 6  14 5 /
\setdashes <1.3mm>
\plot 0 6.5  0 4.3 /
\plot 14 4.5 14 6.5 /

\plot 0 .5  0 3.6 /
\plot 14 .5  14 3.6 /

\setsolid
\Arr{0 4.4}{0 4.3}
\Arr{14 6.4}{14 6.5}

\setsolid
\arr{0.2 1.8}{0.8 1.2}
\arr{0.2 2.2}{0.8 2.8}
\arr{0.2 3.8}{0.8 3.2}
\arr{0.2 4.2}{0.8 4.8}
\arr{0.2 5.8}{0.8 5.2}
\arr{0.2 5}{0.8 5}

\arr{1.2 1.2}{1.8 1.8}
\arr{1.2 2.8}{1.8 2.2}
\arr{1.2 3.2}{1.8 3.8}
\arr{1.2 4.8}{1.8 4.2}
\arr{1.2 5.2}{1.8 5.8}
\arr{1.2 5}{1.8 5}

\arr{2.2 1.8}{2.8 1.2}
\arr{2.2 2.2}{2.8 2.8}
\arr{2.2 3.8}{2.8 3.2}
\arr{2.2 4.2}{2.8 4.8}
\arr{2.2 5.8}{2.8 5.2}
\arr{2.2 5}{2.8 5}

\arr{3.2 1.2}{3.8 1.8}
\arr{3.2 2.8}{3.8 2.2}
\arr{3.2 3.2}{3.8 3.8}
\arr{3.2 4.8}{3.8 4.2}
\arr{3.2 5.2}{3.8 5.8}
\arr{3.2 5}{3.8 5}

\arr{4.2 1.8}{4.8 1.2}
\arr{4.2 2.2}{4.8 2.8}
\arr{4.2 3.8}{4.8 3.2}
\arr{4.2 4.2}{4.8 4.8}
\arr{4.2 5.8}{4.8 5.2}
\arr{4.2 5}{4.8 5}

\arr{5.2 1.2}{5.8 1.8}
\arr{5.2 2.8}{5.8 2.2}
\arr{5.2 3.2}{5.8 3.8}
\arr{5.2 4.8}{5.8 4.2}
\arr{5.2 5.2}{5.8 5.8}
\arr{5.2 5}{5.8 5}

\arr{6.2 1.8}{6.8 1.2}
\arr{6.2 2.2}{6.8 2.8}
\arr{6.2 3.8}{6.8 3.2}
\arr{6.2 4.2}{6.8 4.8}
\arr{6.2 5.8}{6.8 5.2}
\arr{6.2 5}{6.8 5}

\arr{7.2 1.2}{7.8 1.8}
\arr{7.2 2.8}{7.8 2.2}
\arr{7.2 3.2}{7.8 3.8}
\arr{7.2 4.8}{7.8 4.2}
\arr{7.2 5.2}{7.8 5.8}
\arr{7.2 5}{7.8 5}

\arr{8.2 1.8}{8.8 1.2}
\arr{8.2 2.2}{8.8 2.8}
\arr{8.2 3.8}{8.8 3.2}
\arr{8.2 4.2}{8.8 4.8}
\arr{8.2 5.8}{8.8 5.2}
\arr{8.2 5}{8.8 5}

\arr{9.2 1.2}{9.8 1.8}
\arr{9.2 2.8}{9.8 2.2}
\arr{9.2 3.2}{9.8 3.8}
\arr{9.2 4.8}{9.8 4.2}
\arr{9.2 5.2}{9.8 5.8}
\arr{9.2 5}{9.8 5}

\arr{10.2 1.8}{10.8 1.2}
\arr{10.2 2.2}{10.8 2.8}
\arr{10.2 3.8}{10.8 3.2}
\arr{10.2 4.2}{10.8 4.8}
\arr{10.2 5.8}{10.8 5.2}
\arr{10.2 5}{10.8 5}

\arr{11.2 1.2}{11.8 1.8}
\arr{11.2 2.8}{11.8 2.2}
\arr{11.2 3.2}{11.8 3.8}
\arr{11.2 4.8}{11.8 4.2}
\arr{11.2 5.2}{11.8 5.8}
\arr{11.2 5}{11.8 5}

\arr{12.2 1.8}{12.8 1.2}
\arr{12.2 2.2}{12.8 2.8}
\arr{12.2 3.8}{12.8 3.2}
\arr{12.2 4.2}{12.8 4.8}
\arr{12.2 5.8}{12.8 5.2}
\arr{12.2 5}{12.8 5}

\arr{13.2 1.2}{13.8 1.8}
\arr{13.2 2.8}{13.8 2.2}
\arr{13.2 3.2}{13.8 3.8}
\arr{13.2 4.8}{13.8 4.2}
\arr{13.2 5.2}{13.8 5.8}
\arr{13.2 5}{13.8 5}

\endpicture}
$$

Here is the case $n=8$. Now the left hand boundary has to be identified with the right
hand side without any twist. In general, for $n$ even, the Auslander-Reiten quiver of
$\Cal C(H)$ is just $\Bbb ZD_n/\tau^n$.
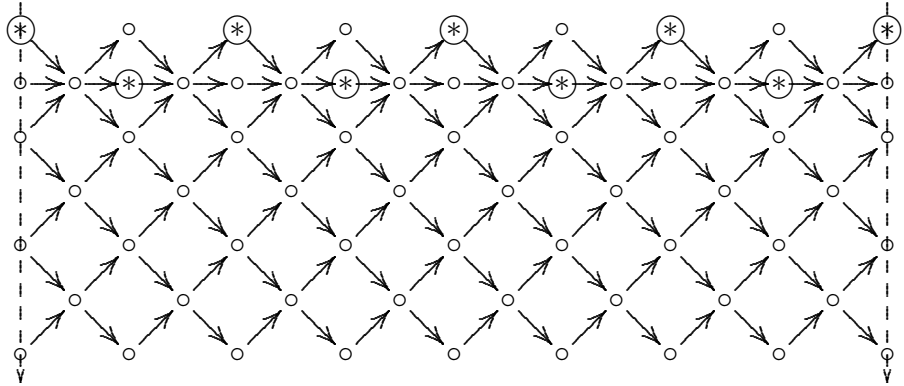
$$
\hbox{\beginpicture
\setcoordinatesystem units <.6cm,.6cm>
\multiput{$\circ$} at 
  0 0  0 2  0 4  0 5  
  1 1  1 3  1 5
  2 0  2 2  2 4       2 6
  3 1  3 3  3 5
  4 0  4 2  4 4  4 5  
  5 1  5 3  5 5
  6 0  6 2  6 4       6 6
  7 1  7 3  7 5
  8 0  8 2  8 4  8 5  
  9 1  9 3  9 5
  10 0  10 2  10 4        10 6
  11 1  11 3  11 5
  12 0  12 2  12 4  12 5  
  13 1  13 3  13 5
  14 0  14 2  14 4        14 6
  15 1  15 3  15 5
  16 0  16 2  16 4  16 5   /
\multiput{$*$} at 0 6  2 5  4 6  6 5  8 6  10 5  12 6  14 5  16 6 /
\multiput{$\bigcirc$} at 0 6  2 5  4 6  6 5  8 6  10 5  12 6  14 5  16 6 /
\setdashes <1.3mm>
\arr{0 6.5}{0 -.5}
\arr{16 6.5}{16 -.5}

\setsolid
\arr{0.2 0.2}{0.8 0.8}
\arr{0.2 1.8}{0.8 1.2}
\arr{0.2 2.2}{0.8 2.8}
\arr{0.2 3.8}{0.8 3.2}
\arr{0.2 4.2}{0.8 4.8}
\arr{0.2 5.8}{0.8 5.2}
\arr{0.2 5}{0.8 5}

\arr{1.2 0.8}{1.8 0.2}
\arr{1.2 1.2}{1.8 1.8}
\arr{1.2 2.8}{1.8 2.2}
\arr{1.2 3.2}{1.8 3.8}
\arr{1.2 4.8}{1.8 4.2}
\arr{1.2 5.2}{1.8 5.8}
\arr{1.2 5}{1.8 5}

\arr{2.2 0.2}{2.8 0.8}
\arr{2.2 1.8}{2.8 1.2}
\arr{2.2 2.2}{2.8 2.8}
\arr{2.2 3.8}{2.8 3.2}
\arr{2.2 4.2}{2.8 4.8}
\arr{2.2 5.8}{2.8 5.2}
\arr{2.2 5}{2.8 5}

\arr{3.2 0.8}{3.8 0.2}
\arr{3.2 1.2}{3.8 1.8}
\arr{3.2 2.8}{3.8 2.2}
\arr{3.2 3.2}{3.8 3.8}
\arr{3.2 4.8}{3.8 4.2}
\arr{3.2 5.2}{3.8 5.8}
\arr{3.2 5}{3.8 5}

\arr{4.2 0.2}{4.8 0.8}
\arr{4.2 1.8}{4.8 1.2}
\arr{4.2 2.2}{4.8 2.8}
\arr{4.2 3.8}{4.8 3.2}
\arr{4.2 4.2}{4.8 4.8}
\arr{4.2 5.8}{4.8 5.2}
\arr{4.2 5}{4.8 5}

\arr{5.2 0.8}{5.8 0.2}
\arr{5.2 1.2}{5.8 1.8}
\arr{5.2 2.8}{5.8 2.2}
\arr{5.2 3.2}{5.8 3.8}
\arr{5.2 4.8}{5.8 4.2}
\arr{5.2 5.2}{5.8 5.8}
\arr{5.2 5}{5.8 5}

\arr{6.2 0.2}{6.8 0.8}
\arr{6.2 1.8}{6.8 1.2}
\arr{6.2 2.2}{6.8 2.8}
\arr{6.2 3.8}{6.8 3.2}
\arr{6.2 4.2}{6.8 4.8}
\arr{6.2 5.8}{6.8 5.2}
\arr{6.2 5}{6.8 5}

\arr{7.2 0.8}{7.8 0.2}
\arr{7.2 1.2}{7.8 1.8}
\arr{7.2 2.8}{7.8 2.2}
\arr{7.2 3.2}{7.8 3.8}
\arr{7.2 4.8}{7.8 4.2}
\arr{7.2 5.2}{7.8 5.8}
\arr{7.2 5}{7.8 5}

\arr{8.2 0.2}{8.8 0.8}
\arr{8.2 1.8}{8.8 1.2}
\arr{8.2 2.2}{8.8 2.8}
\arr{8.2 3.8}{8.8 3.2}
\arr{8.2 4.2}{8.8 4.8}
\arr{8.2 5.8}{8.8 5.2}
\arr{8.2 5}{8.8 5}

\arr{9.2 0.8}{9.8 0.2}
\arr{9.2 1.2}{9.8 1.8}
\arr{9.2 2.8}{9.8 2.2}
\arr{9.2 3.2}{9.8 3.8}
\arr{9.2 4.8}{9.8 4.2}
\arr{9.2 5.2}{9.8 5.8}
\arr{9.2 5}{9.8 5}

\arr{10.2 0.2}{10.8 0.8}
\arr{10.2 1.8}{10.8 1.2}
\arr{10.2 2.2}{10.8 2.8}
\arr{10.2 3.8}{10.8 3.2}
\arr{10.2 4.2}{10.8 4.8}
\arr{10.2 5.8}{10.8 5.2}
\arr{10.2 5}{10.8 5}

\arr{11.2 0.8}{11.8 0.2}
\arr{11.2 1.2}{11.8 1.8}
\arr{11.2 2.8}{11.8 2.2}
\arr{11.2 3.2}{11.8 3.8}
\arr{11.2 4.8}{11.8 4.2}
\arr{11.2 5.2}{11.8 5.8}
\arr{11.2 5}{11.8 5}

\arr{12.2 0.2}{12.8 0.8}
\arr{12.2 1.8}{12.8 1.2}
\arr{12.2 2.2}{12.8 2.8}
\arr{12.2 3.8}{12.8 3.2}
\arr{12.2 4.2}{12.8 4.8}
\arr{12.2 5.8}{12.8 5.2}
\arr{12.2 5}{12.8 5}

\arr{13.2 0.8}{13.8 0.2}
\arr{13.2 1.2}{13.8 1.8}
\arr{13.2 2.8}{13.8 2.2}
\arr{13.2 3.2}{13.8 3.8}
\arr{13.2 4.8}{13.8 4.2}
\arr{13.2 5.2}{13.8 5.8}
\arr{13.2 5}{13.8 5}

\arr{14.2 0.2}{14.8 0.8}
\arr{14.2 1.8}{14.8 1.2}
\arr{14.2 2.2}{14.8 2.8}
\arr{14.2 3.8}{14.8 3.2}
\arr{14.2 4.2}{14.8 4.8}
\arr{14.2 5.8}{14.8 5.2}
\arr{14.2 5}{14.8 5}

\arr{15.2 0.8}{15.8 0.2}
\arr{15.2 1.2}{15.8 1.8}
\arr{15.2 2.8}{15.8 2.2}
\arr{15.2 3.2}{15.8 3.8}
\arr{15.2 4.8}{15.8 4.2}
\arr{15.2 5.2}{15.8 5.8}
\arr{15.2 5}{15.8 5}

\endpicture}
$$

Always, we obtain a serial algebra with $n$ simple modules and the Kupisch series is
$(n-1,\dots,n-1)$. As a quiver with relations, it is given by a cycle with $n$ vertices and
$n$ arrows. If we label all the arrows by $\alpha$, then the relations can be written
just in the form $\alpha^n = 0$. The only difference between the cases $n$ even or odd
concerns the Nakayama permutation: In the even case, it has two orbits, in the odd case only one.
Note that the $A_3$ case may be considered as a $D_3$-case, thus as part of the sequence of serial algebras.

 	\medskip
{\bf 3. The case $D_{2m}$ with $m\ge 3$, using the long arm and one short arm.}
Here we present the case $D_{2m}$ with $m = 4.$
$$
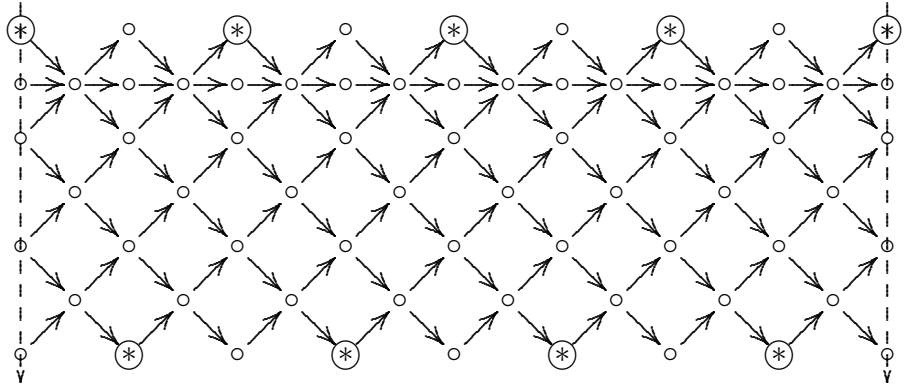
\hbox{\beginpicture
\setcoordinatesystem units <.6cm,.6cm>
\multiput{$\circ$} at 
  0 0  0 2  0 4  0 5  
  1 1  1 3  1 5
       2 2  2 4  2 5  2 6
  3 1  3 3  3 5
  4 0  4 2  4 4  4 5  
  5 1  5 3  5 5
       6 2  6 4  6 5  6 6
  7 1  7 3  7 5
  8 0  8 2  8 4  8 5  
  9 1  9 3  9 5
        10 2  10 4  10 5  10 6
  11 1  11 3  11 5
  12 0  12 2  12 4  12 5  
  13 1  13 3  13 5
        14 2  14 4  14 5  14 6
  15 1  15 3  15 5
  16 0  16 2  16 4  16 5   /
\multiput{$*$} at 0 6    4 6    8 6    12 6  16 6  
  2 0  6 0  10 0  14 0 /
\multiput{$\bigcirc$} at 0 6  4 6   8 6    12 6   16 6 
  2 0  6 0  10 0  14 0 /
\setdashes <1.3mm>
\arr{0 6.5}{0 -.5}
\arr{16 6.5}{16 -.5}

\setsolid
\arr{0.2 0.2}{0.8 0.8}
\arr{0.2 1.8}{0.8 1.2}
\arr{0.2 2.2}{0.8 2.8}
\arr{0.2 3.8}{0.8 3.2}
\arr{0.2 4.2}{0.8 4.8}
\arr{0.2 5.8}{0.8 5.2}
\arr{0.2 5}{0.8 5}

\arr{1.2 0.8}{1.8 0.2}
\arr{1.2 1.2}{1.8 1.8}
\arr{1.2 2.8}{1.8 2.2}
\arr{1.2 3.2}{1.8 3.8}
\arr{1.2 4.8}{1.8 4.2}
\arr{1.2 5.2}{1.8 5.8}
\arr{1.2 5}{1.8 5}

\arr{2.2 0.2}{2.8 0.8}
\arr{2.2 1.8}{2.8 1.2}
\arr{2.2 2.2}{2.8 2.8}
\arr{2.2 3.8}{2.8 3.2}
\arr{2.2 4.2}{2.8 4.8}
\arr{2.2 5.8}{2.8 5.2}
\arr{2.2 5}{2.8 5}

\arr{3.2 0.8}{3.8 0.2}
\arr{3.2 1.2}{3.8 1.8}
\arr{3.2 2.8}{3.8 2.2}
\arr{3.2 3.2}{3.8 3.8}
\arr{3.2 4.8}{3.8 4.2}
\arr{3.2 5.2}{3.8 5.8}
\arr{3.2 5}{3.8 5}

\arr{4.2 0.2}{4.8 0.8}
\arr{4.2 1.8}{4.8 1.2}
\arr{4.2 2.2}{4.8 2.8}
\arr{4.2 3.8}{4.8 3.2}
\arr{4.2 4.2}{4.8 4.8}
\arr{4.2 5.8}{4.8 5.2}
\arr{4.2 5}{4.8 5}

\arr{5.2 0.8}{5.8 0.2}
\arr{5.2 1.2}{5.8 1.8}
\arr{5.2 2.8}{5.8 2.2}
\arr{5.2 3.2}{5.8 3.8}
\arr{5.2 4.8}{5.8 4.2}
\arr{5.2 5.2}{5.8 5.8}
\arr{5.2 5}{5.8 5}

\arr{6.2 0.2}{6.8 0.8}
\arr{6.2 1.8}{6.8 1.2}
\arr{6.2 2.2}{6.8 2.8}
\arr{6.2 3.8}{6.8 3.2}
\arr{6.2 4.2}{6.8 4.8}
\arr{6.2 5.8}{6.8 5.2}
\arr{6.2 5}{6.8 5}

\arr{7.2 0.8}{7.8 0.2}
\arr{7.2 1.2}{7.8 1.8}
\arr{7.2 2.8}{7.8 2.2}
\arr{7.2 3.2}{7.8 3.8}
\arr{7.2 4.8}{7.8 4.2}
\arr{7.2 5.2}{7.8 5.8}
\arr{7.2 5}{7.8 5}

\arr{8.2 0.2}{8.8 0.8}
\arr{8.2 1.8}{8.8 1.2}
\arr{8.2 2.2}{8.8 2.8}
\arr{8.2 3.8}{8.8 3.2}
\arr{8.2 4.2}{8.8 4.8}
\arr{8.2 5.8}{8.8 5.2}
\arr{8.2 5}{8.8 5}

\arr{9.2 0.8}{9.8 0.2}
\arr{9.2 1.2}{9.8 1.8}
\arr{9.2 2.8}{9.8 2.2}
\arr{9.2 3.2}{9.8 3.8}
\arr{9.2 4.8}{9.8 4.2}
\arr{9.2 5.2}{9.8 5.8}
\arr{9.2 5}{9.8 5}

\arr{10.2 0.2}{10.8 0.8}
\arr{10.2 1.8}{10.8 1.2}
\arr{10.2 2.2}{10.8 2.8}
\arr{10.2 3.8}{10.8 3.2}
\arr{10.2 4.2}{10.8 4.8}
\arr{10.2 5.8}{10.8 5.2}
\arr{10.2 5}{10.8 5}

\arr{11.2 0.8}{11.8 0.2}
\arr{11.2 1.2}{11.8 1.8}
\arr{11.2 2.8}{11.8 2.2}
\arr{11.2 3.2}{11.8 3.8}
\arr{11.2 4.8}{11.8 4.2}
\arr{11.2 5.2}{11.8 5.8}
\arr{11.2 5}{11.8 5}

\arr{12.2 0.2}{12.8 0.8}
\arr{12.2 1.8}{12.8 1.2}
\arr{12.2 2.2}{12.8 2.8}
\arr{12.2 3.8}{12.8 3.2}
\arr{12.2 4.2}{12.8 4.8}
\arr{12.2 5.8}{12.8 5.2}
\arr{12.2 5}{12.8 5}

\arr{13.2 0.8}{13.8 0.2}
\arr{13.2 1.2}{13.8 1.8}
\arr{13.2 2.8}{13.8 2.2}
\arr{13.2 3.2}{13.8 3.8}
\arr{13.2 4.8}{13.8 4.2}
\arr{13.2 5.2}{13.8 5.8}
\arr{13.2 5}{13.8 5}

\arr{14.2 0.2}{14.8 0.8}
\arr{14.2 1.8}{14.8 1.2}
\arr{14.2 2.2}{14.8 2.8}
\arr{14.2 3.8}{14.8 3.2}
\arr{14.2 4.2}{14.8 4.8}
\arr{14.2 5.8}{14.8 5.2}
\arr{14.2 5}{14.8 5}

\arr{15.2 0.8}{15.8 0.2}
\arr{15.2 1.2}{15.8 1.8}
\arr{15.2 2.8}{15.8 2.2}
\arr{15.2 3.2}{15.8 3.8}
\arr{15.2 4.8}{15.8 4.2}
\arr{15.2 5.2}{15.8 5.8}
\arr{15.2 5}{15.8 5}

\endpicture}
$$
The quiver of $\widetilde A$ has the following shape (again, the right hand side
has to be identified with the left hand side):
$$
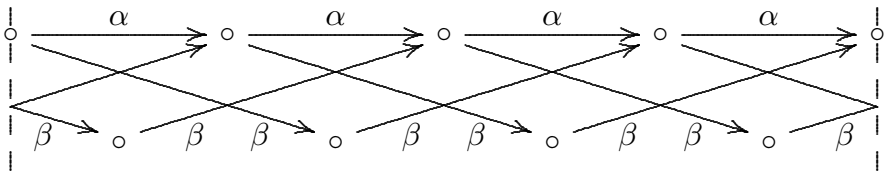
\hbox{\beginpicture
\setcoordinatesystem units <1.2cm,0.6cm>
\multiput{$\circ$} at 0 2  2 2  4 2  6 2  8 2  1 0  3 0  5 0  7 0 /
\arr{0.2 2}{1.8 2}
\arr{2.2 2}{3.8 2}
\arr{4.2 2}{5.8 2}
\arr{6.2 2}{7.8 2}

\arr{0 0.6667}{.8 0.2}
\arr{0.2 1.8}{2.8 0.2}
\arr{2.2 1.8}{4.8 0.2}
\arr{4.2 1.8}{6.8 0.2}
\arr{0 0.667}{1.8 1.8}
\arr{1.2 0.2}{3.8 1.8}
\arr{3.2 0.2}{5.8 1.8}
\arr{5.2 0.2}{7.8 1.8}
\plot 7.2 0.2  8 0.667 /
\plot 6.2 1.8  8 0.667 /
\multiput{$\alpha$} at 1 2.3  3 2.3  5 2.3  7 2.3 /
\multiput{$\beta$} at 0.3 0.1  1.7 0.1
  2.3 0.1  3.7 0.1 
  4.3 0.1  5.7 0.1 
  6.3 0.1  7.7 0.1   /
\setdashes <2mm>
\plot 0 -0.5  0 2.5 /
\plot 8 -0.5  8 2.5 /
\endpicture}
$$
with relations
$$
  \alpha\beta = \beta\alpha = 0,\qquad \alpha^{m-1} = \beta^2.
$$
	\bigskip
\vfill\eject
{\bf References}
	\medskip
\item{[ABS]} Assem, Br\"ustle, Schiffler: Cluster-tilted algebras  
   as trivial extensions.\newline arXiv: math.RT/0601537 

\item{[BMR]} Buan, Marsh, Reiten: Cluster tilted algebras. Trans. Amer. Math. Soc. (To appear). arXiv: math.RT/0402075 

\item{[DR]} Dlab, Ringel: Indecomposable representations of graphs and algebras. Mem. Amer.    Math. Soc. 173 (1976).

\item{[HR]} Happel, Ringel: 
   Tilted algebras. Trans. Amer. Math. Soc. 274 (1982), 399-443. 

\item{[KR]} Keller, Reiten: Cluster-tilted algebras are Gorenstein and stably Calabi-Yau.
   arXiv: math.RT/0512471

\item{[R]} Ringel: Exceptional modules are tree modules. 
    Lin. Alg. Appl. 275-276 (1998) 471-493.
\item{[RV]} Ringel, Vossieck:
   Hammocks. Proc. London Math. Soc. (3) 54 (1987), 216-246. 

\item{[Z]} Bin Zhu: Equivalences between cluster categories. (To appear).
  \newline
   arXiv:math.RT/0511382.

\bye